\documentclass[a4paper,12pt]{amsart}
\usepackage{amsmath,amsthm,amssymb}
\usepackage{mathptmx}
\usepackage{graphicx,graphics}
\usepackage{tikz}
\newtheorem{theorem}{Theorem}[section]
\newtheorem{corollary}{Corollary}[section]
\newtheorem{lemma}{Lemma}[section]

\theoremstyle{definition}
\newtheorem{definition}{Definition}[section]

\newtheorem{example}{Example}[section]

\usepackage{mathrsfs}
\usepackage{enumerate}
\usepackage{float}
\usepackage{graphicx}
\usetikzlibrary{shapes.geometric,arrows}
\usetikzlibrary{mindmap}
\numberwithin{equation}{section}
\begin{document}
\title{A note on distance magic index of partite graphs}

 \author{Eshwar Srinivasan}
 \address{Eshwar Srinivasan, Department of Mathematics, Indian Institute of Technology Guwahati, Guwahati- \textnormal{781039}, India. }
\email{s.eshwar@iitg.ac.in}
\author{A V Prajeesh}
\address{A V Prajeesh, Department of Mathematics, Amrita Vishwa Vidyapeetham 
 Amritapuri-\textnormal{690525}, India.} 
 \email{prajeeshav@am.amrita.edu} 
\author{Krishnan Paramasivam}
 \address{Krishnan Paramasivam, Department of Mathematics, National Institute of Technology Calicut,  Kozhikode-\textnormal{673601}, India. }
 \email{sivam@nitc.ac.in}

\begin{abstract}
	Let $G=(V(G),E(G))$ be a graph on $n$ vertices and let $S$ be a set of any $n$ positive integers. Then $G$ is said to be \textit{$S$-magic} if there exists a bijection $l$ from $V(G)$ to $S$ such that for any vertex $u$ of $G$, $\sum_{v\in N(u)} l(v)$ is a constant $\mu'$ and $\mu'$ is the $S$-magic constant.  Suppose  $\eta(S)$ is the largest positive integer in $S$. If $i(G)$ is the infimum of all $\eta(S)$'s, where the infimum runs overall such sets $S$ for which $G$ admits $S$-magic labeling, then distance magic index $\theta(G)$, of $G$ is defined to be $i(G)-n$. In this article, the distance magic index of certain important classes of partite graphs are determined.
\end{abstract}

\thanks{}

\maketitle
\pagestyle{myheadings}
\markboth{Authors.. }{Distance magic index of graphs}

\section{Introduction}
\noindent In this paper, we consider only simple and finite graphs. We use $V(G)$ for the vertex set and $E(G)$ for the edge set of a graph $G$. The neighborhood, $N_G(v)$ or shortly $N(v)$ of a vertex $v$ of $G$ is the set of all vertices adjacent to $v$ in $G$. For further graph-theoretic terminology and notation, we refer Bondy and Murty \cite{MR2368647} and Hammack $et$ $al.$\cite{MR2817074}.
\par A Kotzig array \cite{kotzig1971magic}, $KA(a,b)$ is an $a\times b$ array in which every row contains each of the integers $0,1,\cdots, b-1$ exactly once and the sum of the entries in each column is the same constant $c = \frac{a(b - 1)}{2}$. The following theorem can be obtained from \cite{MR3013201}.
\par In 2003, Miller {\cite{MR1999203}} introduced the concept of distance magic labeling of a graph $G$, which is injective function $l : V(G)\rightarrow \{1,2,\cdots,|V(G)|\}$ such that for any $u$ of $G$, the weight $w_G(u)$ of $u$ is a constant $\mu$, where $w_G(u)$ is the sum of labels of all neighbors of $u$. A graph $G$ that admits such a labeling $l$ is called a distance magic graph, or shortly, a $dmg$.
\par
\begin{figure}[ht]
	\centering
	\includegraphics[width=120mm]{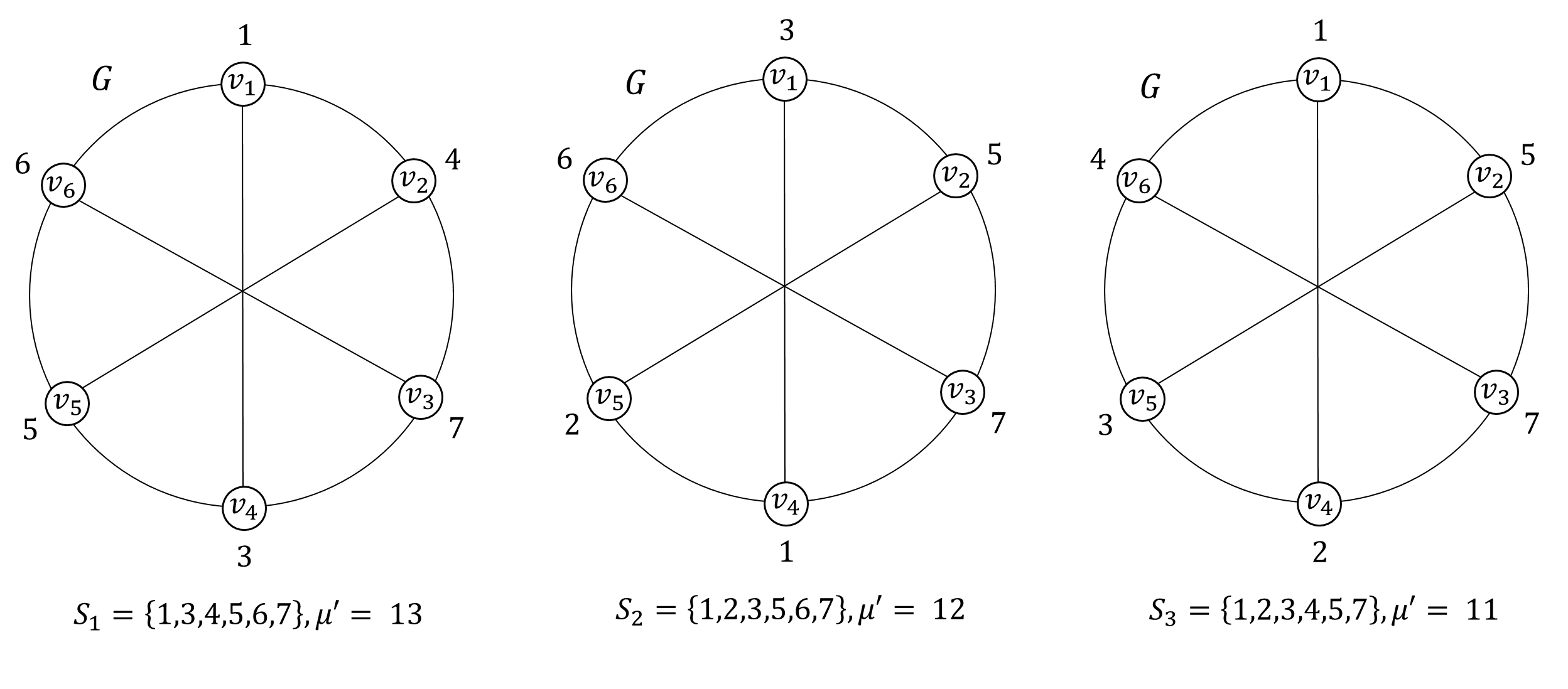}
	\caption{$S_1,$ $S_2$ and $S_3$-magic graph $G$}
	\label{Fig5.1}
\end{figure}

 The following results provide some necessary condition for distance magicness of regular graphs.
\begin{theorem}\label{oddregular}
	\textnormal{\cite{MR1999203,Rao2004,vilfred}} No $r$-regular graph with $r$-odd can be a distance magic graph.
\end{theorem}
\begin{theorem}\label{2mode4}
	\textnormal{\cite{MR2259699}} Let $EIT(a, r)$ be an equalized tournament with an even number $a$ of teams and $r\equiv 2 \mod 4$. Then $a\equiv 0 \mod 4$.
\end{theorem}
\noindent A simple argument on parity, can prove that an odd regular graph is not distance magic \cite{MR1999203} and therefore, $G=K_{3,3}$ given in Figure \ref{Fig5.1}, is not a distance magic. That is, it is not possible to find a labeling from $V(G)$ to $\{1,2,...,6\}$, with all vertices have same constant weight. But the weights of all vertices, are  a unique constant, when the existing set $\{1,2,\cdots,6\}$ is replaced by certain sets $S_1,S_2$ and $S_3$ with $|S_1|=|S_2|=|S_3|=6$ (see Figure 1).
\par Motivated by the above fact, Godinho et al. \cite{MR3743955} defined the concept of $S$-magic labeling of a graph and introduce certain measure namely, the distance magic index to predict the nature of distance of magicness of a graph.  \begin{definition}
	\textnormal{\cite{MR3743955}} Let $G$ be a graph on $n$ vertices and let $S$ be a set of any $n$ positive integers. Then $G$ is said to be \textit{$S$-magic} if there exists a bijection $l$ from $V(G)$ to $S$ such that for any vertex $u$ of $G$, $\sum_{v\in N(u)} l(v)$ is a constant $\mu'$ and $\mu'$ is the $S$-magic constant.
\end{definition}
\begin{definition}
	\textnormal{\cite{MR3743955}} Let $G$ be a graph on $n$ vertices and let $S$ be a set of $n$ positive integers with the largest number $\eta(S)$. Define $i(G)$ to be the infimum of all $\eta(S)$'s, where the infimum runs overall $S$ for which $G$ admits $S$-magic labeling. Then the \textit{distance magic index} of $G$, $\theta(G)=i(G)-n$.
\end{definition}
In 2015, Godinho et al. \cite{godinho2015s} defined the concept of  fractional domination number and its relationship with $S$-magic graphs.
\begin{definition}\cite{godinho2015s}
    Let $G = (V,E)$ be a graph without isolated vertices. A function $g : V \rightarrow [0, 1]$ is called a total dominating function of $G$ if $g(N(v)) = \sum_{u \in N(v)}  g(u) \ge 1$ for all $v \in V$. The fractional total domination number $\gamma_{ft}(G)$ is defined to be $min\{|g| : g$ is a total dominating function of $G\}$, where $|g| = \sum_{v \in V} g(v)$.
    \end{definition}
 
\begin{theorem}\textnormal{{\cite{MR1999203,kotlar}}}     Let $1 \le n_1 \le n_2 \le n_3 \le n_4$ and $s_{i} = \sum_{i = 1}^{4} n_i$. Suppose $G$ is either $K_{n_1,n_2}$ or $K_{n_1,n_2,n_3}$ or $K_{n_1,n_2,n_3, n_4}$. There exists a distance magic labeling of $G$ if and only if the following conditions hold.
    \begin{enumerate}
    \item[$(a)$] $n_2 \ge 2$,
    \item[$(b)$] $n(n + 1) \equiv 0$ (mod $2p)$, where $n = s_{p} = |V (K_{n_{1},,...,n_{p}})|$, and
    \item[$(c)$] $\displaystyle \sum_{j=1}^{s_{i}} (n - j + 1) \ge \frac{n(n+1)i}{2p}$ for $1 \le i \le p$,
    \end{enumerate}
    where $p=2$ if $G$ is complete bipartite graph, and $p=3$ if $G$ is complete tripartite graph, and $p=4$ if $G$ is complete four-partite graph.
 \end{theorem} 
 \begin{lemma} \textnormal{ \cite{godinho2015s})}
     The complete r-partite graph $G = K_{m_1,m_2,...,n_r}$ is $S$-magic with labeling f if and only if the sum of the labels of all the vertices in any two partite sets are equal.
 \end{lemma}
 \begin{theorem} \textnormal{\cite{MR3743955}}
Let $G$ be the complete bipartite graph $K_{n_{1}, n_{2}}$, where $2 \le n_{1} \le n_{2}$ and $n_{1} + n_{2} = n$. Then\\
    $\displaystyle \theta(G)$ = $\begin{cases}
    0 & \text{ if } n(n+1) \ge 2n_{2}(n_{2} + 1) \textnormal{\it ~and~} n\equiv 0 \text{ or } 3 (\text{mod~} 4)\\
    1 & \text{ if } n(n+1) \ge 2n_{2}(n_{2} + 1) \textnormal{\it ~and~} n\equiv 1 \text{ or } 2 (\text{mod~} 4)\\
     \bigl\lceil \frac{|n(n+1) - 2n_{2}(n_{2} + 1)|}{2n_
     {1}} \bigr\rceil & \text{ if } n(n+1) < 2n_{2}(n_{2} + 1).
    \end{cases}$
 \end{theorem}
    \begin{theorem} \textnormal{\cite{godinho2015s}}
       If a graph $G$ admits an $S$-magic labeling $f$ with magic constant $k$, then $k = \frac{\alpha}{\gamma_{ft}(G)}$, where $\alpha = \sum_{i \in S}i$.
    \end{theorem}
    \begin{corollary} \textnormal{\cite{godinho2015s}}
    If $G$ is $S$-magic, then the smallest $S$-magic constant corresponds to the $S$-magic labeling for which $  \sum_{i \in S}i$ is minimum.
    \end{corollary}
    \par Let $G$ be a graph on $n$ vertices. It is clear that $G$ is distance magic if and only if $\theta(G) = 0$. On the other hand, if $G$ is not $S$-magic, for any $S$ with $|V(G)|$ positive integers, then $\theta(G)$ is not finite. Godinho et al.\cite{MR3743955} identified a parameter for a general graph $G$, defined by $g(x) = \frac{1}{2}[(2\delta(n+x)-\delta^2+\delta)-\Delta(\Delta+1)]$, which can be effectively used to determine the distance magic index of general graphs. From \cite{MR3743955}, one can see that if $g(0)<0$, then $G$ is not a distance magic graph and hence $\theta(G)\geq 1$. 
The following result gives a lower bound for $\theta(G)$.
\begin{theorem}\label{godinho}
	\textnormal{\cite{MR3743955}} If $G$ is a graph of order $n$ such that $g(0)<0$, then $\theta(G)\geq \biggl\lceil \frac{|g(0)|}{\delta}\biggr\rceil$.
\end{theorem}
\par However, if $g(0)\geq 0$, there is no result on $G$ to argue whether for that graph, $\theta(G)=0$ or $\theta(G) \geq 1$. Note that if $G$ is an $r$-regular graph with $n$ vertices, then $g(0)>0$ and hence Theorem \ref{godinho} does predict any possible $\theta(G)$. 
\section{Distance magic index of complete-tripartite graphs}
In this section, we prove the distance magicness and the magic index of the complete tripartite graph $K_{n_1, n_2, n_3}$, where $V=\bigcup_{i=1}^3 V_i$ and $|V_i|=n_i$ so that $n_1 + n_2 + n_3 = |V|=n$ (say) and the edge set $E =\{uv: u\in V_i,v\in V_j, 1\leq i,j\leq 3 \textnormal{~and~} i\ne j\}$ with $|E |=n_1n_2 + n_2n_3 + n_1n_3$.    Also when $l$ is a mapping from the vertex set of $K_{n_1, n_2, n_3}$ to $S$, where $S$ is set of positive integers with $|S|=|V(K_{n_1, n_2, n_3})|$, we say label set of $V_i$ is $\{l(v) : v\in V_i\}$, which is denoted by $L_i$ and  we use $s(L_i)$ to denote the sum of labels of all the elements in $L_i$, where $i=1,2,3$. Without loss of generality, one can restrict $2 \le n_{1} \le n_{2} \le n_{3}$. For any $i<j$, we denote $\zeta_i^j = \sum_{\gamma=i}^{j} \gamma$. Now, by using the above notation and terminology, we first rewrite the third condition in Theorem 1.3\cite{MR1999203, kotlar} for the complete tripartite graph $K_{n_1, n_2, n_3}$
 \begin{equation}\label{1a}
     \zeta_{n - n_1 + 1}^{n} = \sum_{j=1}^{n_{1}} (n-j+1) = \frac{n_{1}}{2} (2n - n_{1} + 1) \ge \frac{n(n+1)}{6} = \frac{\zeta_1^{n}}{3}  
 .\end{equation}
 and, 
 \begin{equation*}
      \zeta_{n_3 + 1}^{n} = \sum_{j=1}^{n_{1}+n_{2}} (n-j+1) = \frac{n_{1}+n_{2}}{2} (2n - n_{1} - n_{2} + 1) \ge \frac{n(n+1)}{3} = \frac{2\zeta_1^{n}}{3}
 .\end{equation*}
 Since,
 \begin{equation}\label{2a}
     \frac{n_{1}+n_{2}}{2} (2n - n_{1} - n_{2} + 1) + \frac{n_{3}(n_{3}+1)}{2} = \frac{\zeta_1^{n}}{3} + \frac{2\zeta_1^{n}}{3} = \frac{n(n+1)}{2} = \zeta_1^{n} 
     .\end{equation}
     we get,
    \begin{equation}\label{3a}
        3\zeta_1^{n_3} \le \zeta_1^{n}
    .\end{equation}
    From the above equations \eqref{1a} and \eqref{3a}, one can write the following theorem.
    \begin{theorem}
        $K_{n_{1}, n_{2}, n_{3}}$ is distance magic if and only if the following conditions hold.
    \begin{enumerate}
    \item[$1)$] $n_{2} \ge 2$,
    \item[$2)$] $2\zeta_1^n \equiv 0\mod 6$,
    \item[$3)$] $3\zeta_{n-n_1+1}^{n} \ge \zeta_1^{n}$ and,
    \item[$4)$] $3\zeta_{1}^{n_3} \le \zeta_1^{n}$.
    \end{enumerate}
    \end{theorem}
    Note that  Theorem [3.1]  is a rewritten from Theorem [1.3] for the particular case of complete tripartite graphs, $K_{n_{1}, n_{2}, n_{3}}$. Therefore, from above theorem, one can classify $K_{n_{1}, n_{2}, n_{3}}$ into the following five cases.
    \begin{enumerate}
        \item[$i)$] $3\zeta_{n-n_1+1}^{n} \ge  \zeta_1^{n} \ge 3\zeta_{1}^{n_3}   \text{ and }2\zeta_1^{n}  \equiv 0 \text{ mod } 6$.
        \item[$ii)$] $\zeta_1^{n} > 3\zeta_{n-n_1+1}^{n}$ and $\zeta_1^{n} \ge 3\zeta_{1}^{n_3} $.
        \item[$iii)$] $3\zeta_{n-n_1+1}^{n} < \zeta_1^{n}  < 3\zeta_{1}^{n_3}   $.
        \item[$iv)$]  $3\zeta_{n-n_1+1}^{n} \ge  \zeta_1^{n} \ge 3\zeta_{1}^{n_3}  \text{ and } 2\zeta_1^{n}  \equiv 2 \text{ mod } 6$.
        \item[$v)$]  $\zeta_1^{n} \le 3\zeta_{n-n_1+1}^{n}$ and $\zeta_1^{n} < 3\zeta_{1}^{n_3} $.
    \end{enumerate}
    
    The resulting theorem determines the distance magic index of the complete tripartite graphs for some of the cases mentioned above.
    \begin{lemma}
        Let $K_{n_{1}, n_{2}, n_3}$ be a complete tripartite graph such that $$3\zeta_{n-n_1+1}^{n} \ge  \zeta_1^{n} \ge 3\zeta_{1}^{n_3}  \text{ and } 2\zeta_1^{n} \equiv 2 \mod 6.$$ Then $K_{n_1 + 1, n_2, n_3}$ is a distance magic graph and consequently, $\theta(K_{n_1 + 1, n_2, n_3}) = 0$.
    \end{lemma}
    \begin{proof}    Consider, $6\zeta_{n-n_1+1}^{n+1} = 3(n_1 + 1)(2(n+1) - (n_1+1) + 1)$,
       \begin{eqnarray*}
           3(n_1 + 1)(2(n+1) - (n_1+1) + 1) &=& 3n_1(2n - n_1 + 1) - 6n_1 + 6n\\ &\ge& n(n+1) - 6n_1 + 6n.
       \end{eqnarray*}
        Since $n_1 \le \frac{n}{3}$. Hence,
        \begin{eqnarray*}
            3(n_1 + 1)(2(n+1) - (n_1+1) + 1) &\ge& n(n+1) + 6n - 6n_1\\ &\ge& n(n+1) + 4n\\ &=& n^2 + 5n \\ &>& n^2 + 3n + 2\\ &=& (n+1)(n+2) = 2\zeta_{1}^{n+1}.
        \end{eqnarray*}
        Now consider $6\zeta_{1}^{n_3} = 3n_3(n_3+1)$,
        \begin{eqnarray*}
            3n_3(n_3 + 1) \le n(n+1) \le (n+1)(n+2) = 2\zeta_{1}^{n+1}.
        \end{eqnarray*}
    
       \noindent Therefore, 
        \begin{eqnarray*}
            2\zeta_{1}^{n+1} = (n+1)(n+2)  \cong n(n+1) + 2(n+1) \text{ mod } 6 \cong 2n + 4 \text{ mod } 6.
        \end{eqnarray*}
        Since $2\zeta_{1}^{n} \cong 2$ mod $6$, $2n \cong 2$ mod $6$, we have,
        \begin{equation*}
            2\zeta_{1}^{n+1} \cong 0 \text{ mod } 6
        .\end{equation*}
        Hence by Theorem 3.1, $K_{n_1 + 1, n_2, n_3}$ is a distance magic graph.
        \end{proof}
    \begin{theorem}
       Let G be the complete tripartite graph $K_{n_{1}, n_{2}, n_{3}}$. Then,
\begin{enumerate}
    \item[$(i)$] If  $3\zeta_{n-n_1+1}^{n} \ge  \zeta_1^{n} \ge 3\zeta_{1}^{n_3}   \text{ and }2\zeta_1^{n}  \equiv 0 \text{ mod } 6$, then $\theta(G)=0$. 
       
    \item[$(ii)$] If $\zeta_1^{n} > 3\zeta_{n-n_1+1}^{n}$ and $\zeta_1^{n} \ge 3\zeta_{1}^{n_3} $, then
    \begin{equation*}
        \theta(G) = \begin{cases}
         \biggl\lceil \frac{3\zeta_{1}^{n-n_1} - 2\zeta_{1}^{n}}{2n_{1}} \biggr\rceil & \text{ if } n - n_{1} \equiv 0 \text{ or } 3 \text{ (mod } 4)\\
        \biggl\lceil \frac{3\zeta_{1}^{n-n_1} - 2\zeta_{1}^{n} + 1}{2n_{1}}\biggr \rceil & \text{ if } n - n_{1} \equiv 1 \text{ or } 2 \text{ (mod } 4). \end{cases} .\end{equation*}
    \item[$(iii)$] If
   $3\zeta_{n-n_1+1}^{n} < \zeta_1^{n}  < 3\zeta_{1}^{n_3}$, then
    \begin{equation*}
        \theta(G) \ge  \biggl\lceil \frac{\zeta_{1}^{n_3} - \zeta_{n-n_1+1}^{n}}{n - n_{3}} \biggr\rceil.
    \end{equation*}
    \item[$(iv)$] If $3\zeta_{n-n_1+1}^{n} \ge  \zeta_1^{n} \ge 3\zeta_{1}^{n_3}  \text{ and } 2\zeta_1^{n}  \equiv 2 \text{ mod } 6$,
    \begin{equation*}
        1 \le \theta(G) \le n+1.
    \end{equation*}
    \item[$(v)$] If $\zeta_1^{n} \le 3\zeta_{n-n_1+1}^{n}$ and $\zeta_1^{n} < 3\zeta_{1}^{n_3},$ then
    \begin{equation*}
        \theta(G) =  \biggl\lceil \frac{\zeta_{1}^{n_3} - \zeta_{n-n_1+1}^{n}}{n_{1}} \biggr\rceil,  \text{ if } \frac{\zeta_{1}^{n_3} - \zeta_{n-n_1+1}^{n}}{n_{1}} \ge \theta(H)
    \end{equation*}
    \begin{equation*}
        \theta(G) \le \theta(H) + 1,  \text{ if } \frac{\zeta_{1}^{n_3} - \zeta_{n-n_1+1}^{n}}{n_{1}} < \theta(H),\textnormal{~~~~~~~~}
    \end{equation*}\\
    where $H \cong K_{n_{2}, n_{3}}$.
    
\end{enumerate}
 \end{theorem}
 \begin{proof}
Let $G=K_{n_{1}, n_{2}, n_{3}}$ be a complete tripartite graph, where  $2 \le n_{1} \le n_{2}\le n_{3}$ and $V =V_1\cup V_2\cup V_3$ with $|V_i|=n_i$, and  $|V | =n_1+n_2+n_3=n$~(say) and with the edge set $E =\{uv: u\in V_i,v\in V_j, 1\leq i,j\leq 3 \textnormal{~and~} i\ne j\}$. Note that $|E|=n_1n_2 + n_2n_3 + n_1n_3$.  For any vertex $v\in V_i$, since the neighborhood of $v$ is $V_j\cup V_k$, the degree of $v$ is  $n_j+n_k$, where $i,j,k$ are three distinct integers and $1 \leq i,j,k\leq 3$.\\
    \noindent \textbf{Proof of case $(i)$} \\ Suppose that $3\zeta_{n-n_1+1}^{n} \ge  \zeta_1^{n} \ge 3\zeta_{1}^{n_3} \text{ and}2\zeta_1^{n}  \equiv 0 \mod 6$. By using Theorem 3.1, we have $\theta(G) = 0$.
    
     \noindent \textbf{Proof of case $(ii)$} \\ Let $L_{i}$ denote the set of labels of the partition $V_{i}$, where $i=1,2,3$.     
    Let the label sets of partitions $V_{1}$, $V_{2}$ and $V_{3}$ be $L_{1} = \{n - n_{1} + 1, \cdots, n\}$, $L_{2} = \{n_{3} + 1,\cdots, n - n_{1}\}$ and $L_{3} = \{1,\cdots, n_{3}\}$ respectively. Let $H$ be a subgraph of $G$ such that $H \cong K_{n_{2},n_{3}}$. Consider,
    \begin{eqnarray*}
        \frac{(n - n_{1})(n - n_{1} + 1)}{2} = \frac{n^{2} + n + n_{1}^{2} - 2nn_{1} - n_{1}}{2}\\
        = \frac{n^{2} + n}{2} - \frac{n_{1}(2n - n_{1} + 1)}{2}.
    \end{eqnarray*}
    Since $3n_{1}(2n - n_{1} +1) < n(n+1)$ and $3n_{3}(n_{3} + 1) \le n(n+1)$,
    \begin{eqnarray}\label{6}
        \frac{(n - n_{1})(n - n_{1} + 1)}{2} > \frac{n^{2} + n}{3} \ge n_{3}(n_{3} + 1).
    \end{eqnarray}
    Hence by Theorem 1.4,\\
    \begin{equation}\label{7}
        \theta(H) = \begin{cases}
    0 & \text{if}  n - n_{1} \equiv 0 \text{or} 3 (\mod 4)\\
    1 & \text{ if }  n - n_{1} \equiv 1 \text{ or } 2 (\mod 4).\\
    \end{cases} \\
    \end{equation}
    \vskip .2cm \noindent Hence by Theorem 1.4 \cite{MR3743955}, the graph $H \equiv K_{n_{2},n_{3}}$ is $S$-magic with $S = \{1, 2,\cdots, n-n_{1}\}$ if $\theta(H) = 0$ and with $S = \{1, 2,\cdots, n-n_{1}-1, n-n_{1}+1\}$ if $\theta(H) = 1$. Hence from this,
\begin{eqnarray}\label{8}
       h = \begin{cases} \frac{(n-n_{1})(n-n_{1}+1)}{4}, & \textnormal{if} n-n_{1} \equiv 0 \textnormal{ or} 3 (\mod 4)\\
        \frac{(n-n_{1})(n-n_{1}+1)+2}{4}, & \textnormal{if} n-n_{1} \equiv 1 \text{or} 2 (\mod 4).
        \end{cases}
\end{eqnarray}
    where $h$ is the magic constant of the subgraph $H \cong K_{n_{2}, n_{3}}$. Note that $h$ is an integer.\newline
    
    \noindent \textbf{Subcase 1} Let $n - n_{1} \equiv 0$ or $3(\mod 4)$. Then from \eqref{7}, $\theta(H) = 0$. Thus $L_{1} = \{n - n_{1} + 1,\cdots, n\}$, $L_{2}$ is the set of labels of  $V_{2}$ in $H$ and $L_{3}$ is the set of all labels of $V_{3}$ in $H$. Now for the graph $G$ to be $S$-magic, $s(L_{1}) = h$. But $s(L_{1}) = \frac{n_{1}}{2}(2n - n_{1} + 1)$. Therefore, let
    \begin{eqnarray*}
    \centering
     \lambda &=& h - s(L_1) = \frac{(n-n_1)(n-n_1+1)}{4} - \frac{n_1}{2}(2n - n_1 + 1).\\
    \textnormal{Then, }\mbox{~~~~~~~~~} \lambda &=& h - \{s(L_1) + s(L_2) + s(L_3)\} + s(L_2) + s(L_3)\\ 
          &=& 3h - \frac{n(n+1)}{2}.
    \end{eqnarray*}
    \begin{equation}\label{9}
        \therefore \lambda = \frac{6h - n(n+1)}{2} =  \frac{3(n-n_{1})(n-n_{1}+1)-2n(n+1)}{4}
    .\end{equation}
    From equation \eqref{6}, $\lambda > 0$. Hence, let $\lambda = n_{1}q + r$, where $q \ge 0$ and $0 \le r < n_{1}$. Now, increase each label of $L_{1}$ by $q$ and keep $L_{2}$ and $L_{3}$. Then, $L_{1}^{'} = \{n - n_{1} + q + 1,\cdots, n + q\}$, $L_{2}^{'} = L_{2}$ and $L_{3}^{'} = L_{3}$. 
    \vskip .2cm \noindent Therefore, $$s(L_{1}^{'}) = s(L_{1}) + n_{1}q.$$
   \noindent {\bf (i)} Suppose $r=0$. We can label the partitions $V_{1}$, $V_{2}$ and $V_{3}$ with the labels $L_{1}^{'}$, $L_{2}^{'}$ and $L_{3}^{'}$ respectively. From equation \eqref{8} and \eqref{9}, it is clear that, $s(L_{1}^{'}) = s(L_{2}^{'}) = s(L_{3}^{'})$. Hence by Lemma 1.1, $G$ is $S$-magic. Since the highest label is $n+q$, the magic index, $\theta(G) = q$.\\
    {\bf (ii)} Suppose $r>0$, We replace the label $n + q - r + 1 \in L_{1}^{'}$ with the label $n + q - r + 1 + r = n + q + 1$. Hence the resulting labels are $L_{1}^{''} = \{n - n_{1} + q + 1,\cdots, n + q - r, n + q - r +2,\cdots, n + q + 1\}$, $L_{2}^{''} = L_{2}^{'}$ and $L_{3}^{''} = L_{3}^{'}$. Now $s(L_{1}^{''}) = s(L_{1}^{'}) + r = s(L_{1}) + n_{1}q + r$, $s(L_{2}^{''})=s(L_{2})$ and $s(L_{3}^{''}) = s(L_{3})$. We label the partitions $V_{1}$, $V_{2}$ and $V_{3}$ with the labels $L_{1}^{''}$, $L_{2}^{''}$ and $L_{3}^{''}$ respectively. 
    \vskip .2cm \noindent From equation \eqref{8} and \eqref{9}, it is clear that $s(L_{1}^{''}) = s(L_{2}^{''}) =s(L_{3}^{''})$. Hence by Lemma 1.1, $G$ is $S$-magic. Since the highest label is $n + q + 1$, the magic index, $\theta(G) = q + 1$. Therefore,
    \begin{equation*}
        \theta(G) = \biggl\lceil \frac{6h - n(n+1)}{2n_{1}} \biggr\rceil = \biggl\lceil \frac{3(n - n_{1})(n - n_{1} + 1) - 2n(n+1)}{4n_{1}} \biggr\rceil
    .\end{equation*}
    Since the subgraph $H$ is $S$-magic with $\theta(H) = 0$, by Corollary 1.1, $h$ is the minimum magic constant for the subgraph $H$. Therefore, $\theta(G)$ is minimum as it is a linear function of $h$. \newline
    \noindent \textbf{Subcase 2}. Let $n - n_{1} \equiv 1 $ or $2$ (mod $4)$. Hence from \eqref{7}, $\theta(H) = 1$. So, take $L_{1}^{'} = \{n-n_{1}, n-n_{1}+2,\cdots, n\}$, $L_{2}^{'}$ is the set of labels of $V_{2}$ in $H$ and $L_{3}^{'}$ is the set of labels of $V_{3}$ in $H$. Now if $G$ is $S$-magic, then by Lemma 1.1, we must have $s(L_{1}^{'}) = s(L_{2}^{'}) = s(L_{3}^{'})$. But $s(L_{1}^{'}) = s(L_{1}) - 1 = \frac{n_{1}}{2}(2n - n_{1} +1) - 1$. Hence, let
    \begin{eqnarray*}
        \lambda &=& h - s(L_{1}) + 1 = \frac{(n-n_{1})(n-n_{1}+1)+2}{4} - \bigl(\frac{n_{1}}{2}(2n - n_{1} + 1) -1\bigr)   \\
        \implies \lambda &=& h - s(L_{1}) + 1 - s(L_{2}) - s(L_{3}) + s(L_{2}) + s(L_{3}) \\ &=& 3h - \frac{n(n+1)}{2}.\\  
       \textnormal{Thus,~~} \lambda &=& \frac{6h - n(n+1)}{2} =  \frac{3(n-n_{1})(n-n_{1}+1)-2n(n+1) + 6}{4}
    \end{eqnarray*}
    \\ \noindent  As in the previous subcase, from equation \eqref{6}, $\lambda > 0$. Hence, let $\lambda - 1 = n_{1}q + r - 1$, where $q \ge 0$ and $0 \le r < n_{1}$. 
    \vskip .2cm \noindent  First, let us consider $q > 1$. Now  increase each label of $L_{1}^{'}$ by $q$ and keep $L_{2}^{'}$ and $L_{3}^{'}$ unchanged, we get $L_{1}^{''} = \{n - n_{1} + q, n - n_{1} + 2 + q,\cdots, n + q\}$, $L_{2}^{''} = L_{2}$ and $L_{3}^{''} = L_{3}$. Therefore, $s(L_{1}^{''}) = s(L_{1}) + n_{1}q - 1$. If $r = 0$, we can label the partitions $V_{1}$, $V_{2}$ and $V_{3}$ with the labels $L_{1}^{''}$, $L_{2}^{''}$ and $L_{3}^{''}$ respectively. From equation \eqref{8} and \eqref{9}, we see that $s(L_{1}^{''}) = s(L_{2}^{''}) = s(L_{3}^{''})$. Hence by Lemma 1.1, $G$ is $S$-magic. Since the highest label is $n + q$, the magic index, $\theta(G) = q$. Since $\lambda - 1 = n_{1}(q - 1) + n_{1} - 1$, we have the distance magic index, $$\theta(G) = q = \biggl\lceil \frac{\lambda - 1}{n_{1}} \biggr\rceil.$$
    
    \vskip .2cm \noindent {\bf (i)} Suppose $r$ is $1$. Then replace the label $n - n_{1} + q \in L_{1}^{''}$ with the label $n - n_{1} + q +1$. Hence the resulting labels are $L_{1}^{'''} = \{n - n_{1} + q + 1, n - n_{1} + q + 2, \dots, n + q\}$, $L_{2}^{'''} = L_{2}$ and $L_{3}^{'''} = L_{3}$. Therefore, $s(L_{1}^{'''}) = s(L_{1}) + n_{1}q$, $s(L_{2}^{'''}) = s(L_{2})$ and $s(L_{3}^{'''}) = s(L_{3})$. We label the partitions $V_{1}$, $V_{2}$ and $V_{3}$ with the labels $L_{1}^{'''}$, $L_{2}^{'''}$ and $L_{3}^{'''}$ respectively. From equation \eqref{8} and \eqref{9}, we see that $s(L_{1}^{'''}) = s(L_{2}^{'''}) = s(L_{3}^{'''})$. Hence, by Lemma 1.1, $G$ is $S$-magic. Since the highest label is $n + q$, the distance magic index $$\theta(G) = q = \biggl\lceil \frac{\lambda - 1}{n_{1}} \biggr\rceil.$$
    
     \vskip .2cm \noindent {\bf (ii)} Suppose $r>1$.  Then replace the label $n - n_{1} + q \in L_{1}^{''}$ with the label $n - n_{1} + q +1$ and $n + q - r + 2 \in L_{1}^{''}$ with the label $n + q - r + 2 + r - 1 = n + q + 1$. Hence the resulting labels are $L_{1}^{''''} = \{n - n_{1} + 1 + q, n - n_{1} + 2 + q, \dots, n + q - r + 1, n + q - r + 3, ..., n + q + 1\}$, $L_{2}^{''''} = L_{2}^{'}$ and $L_{3}^{''''} = L_{3}^{'}$. Now $s(L_{1}^{''''}) = s(L_{1}^{'}) + r = s(L_{1}) + n_{1}q + r - 1$, $s(L_{2}^{''''}) = s(L_{2})$ and $s(L_{3}^{''''}) = s(L_{3})$. We label the partitions $V_{1}$, $V_{2}$ and $V_{3}$ with the labels sets $L_{1}^{''''}$, $L_{2}^{''''}$ and $L_{3}^{''''}$ respectively. From equation \eqref{8} and \eqref{9}, we see that  $s(L_{1}^{''''}) = s(L_{2}^{''''}) = s(L_{3}^{''''})$. Hence, by Lemma 1.1, $G$ is $S$-magic. Since the highest label is $n + q + 1$, the distance magic index $\theta(G) = q + 1$. Since $\lambda - 1 = n_{1}q + r - 1$ and $r > 1$, $$\theta(G) = q + 1 = \biggl\lceil \frac{\lambda - 1}{n_{1}} \biggr\rceil.$$
    
    \noindent Now, if we take $q = 1$, then $\lambda  - 1 = n_{1} + r - 1$. Hence, we can not increase  $n - n_{1} \in L_{1}^{'}$ by 1 as $n - n_{1} + 1 \in L_{2}^{'} \cup L_{3}^{'}$. Hence, increase each label of $L_{1}^{'}-\{n-n_1\}$ by 1, and keep $L_{2}^{'}$ and $L_{3}^{'}$ unchanged and hence, we get, $L_{1}^{'''} = \{n - n_{1}, n - n_{1} + 3,\cdots, n + 1 \}$, $L_{2}^{'''} = L_{2}$ and $L_{3}^{'''} = L_{3}$. Therefore $$s(L_{1}^{'''}) = s(L_{1}) + n_{1} - 1.$$ 
     \vskip .2cm \noindent {\bf (i)} Suppose $r=0$. Then replace the label $n + 1 \in L_{1}^{'''}$ with $n + 2$. Hence the resulting labels are $L_{1}^{''''} = \{n - n_{1}, n - n_{1} + 3,\cdots, n, n + 2\}$, $L_{2}^{''''} = L_{2}$ and $L_{3}^{''''} = L_{3}$. Therefore, $s(L_{1}^{''''}) = s(L_{1}) + n_{1}$, $s(L_{2}^{''''}) = s(L_{2})$ and $s(L_{3}^{''''}) = s(L_{3})$. We label the partitions $V_{1}$, $V_{2}$ and $V_{3}$ with the labels $L_{1}^{''''}$, $L_{2}^{''''}$ and $L_{3}^{''''}$, respectively. From equation \eqref{8} and  equation \eqref{9}, it is clear that, $s(L_{1}^{''''}) = s(L_{2}^{''''}) = s(L_{3}^{''''})$. Hence by Lemma 1.1, $G$ is $S$-magic. Since the highest label is $n + 2 = n + q + 1$, the magic index, $\theta(G) = 2 = q+1$. Since $\lambda - 1 = n_{1} + r - 1$, we have $\theta(G) = \lceil \frac{\lambda - 1}{n_{1}} \rceil$.\\
     \vskip .2cm \noindent {\bf (ii)} Suppose $r = 1$. Then replace the label $n - n_{1} \in L_{1}^{'''}$ with $n - n_{1} + 2$. Hence the resulting labels are $L_{1}^{''''} = \{n - n_{1} + 2, n - n_{1} + 3,\cdots, n + 1\}$, $L_{2}^{''''} = L_{2}$ and $L_{3}^{''''} = L_{3}$. Therefore, $s(L_{1}^{''''}) = s(L_{1}) + n_{1}$, $s(L_{2}^{''''}) = s(L_{2})$ and $s(L_{3}^{''''}) = s(L_{3})$. We label the partitions $V_{1}$, $V_{2}$ and $V_{3}$ with the label set $L_{1}^{''''}$, $L_{2}^{''''}$ and $L_{3}^{''''}$, respectively. From \eqref{8} and \eqref{9} it is clear that, $s(L_{1}^{''''}) = s(L_{2}^{''''}) = s(L_{3}^{''''})$. Hence by lemma 1.1, $G$ is $S$-magic. Since the highest label is $n + 1 = n + q $, the magic index, $\theta(G) = 1 = q$. Since $\lambda - 1 = n_{1}$, $\theta(G) = \lceil \frac{\lambda - 1}{n_{1}} \rceil$.\\
     \vskip .2cm \noindent {\bf (iii)} Suppose that $r > 1$.  We replace the label $n - n_{1} \in L_{1}^{'''}$ with the label $n - n_{1} + 2$ and $n - r + 3 \in L_{1}^{'''}$ with $n - r + 3 + r - 1 = n + 2$. Hence the resulting label sets are $L_{1}^{''''} = \{n - n_{1} + 2, n - n_{1} + 3,\cdots, n - r + 2, n - r + 4,\cdots, n +2\}$, $L_{2}^{''''} = L_{2}$ and $L_{3}^{''''} = L_{3}$. Now $s(L_{1}^{''''}) = s(L_{1}) + n_{1} + r$, $s(L_{2}^{''''}) = s(L_{2})$ and $s(L_{3}^{''''}) = s(L_{3})$. We label the partitions $V_{1}$, $V_{2}$ and $V_{3}$ with the labels $L_{1}^{''''}$, $L_{2}^{''''}$ and $L_{3}^{''''}$ respectively. From \eqref{8} and \eqref{9} it is clear that, $s(L_{1}^{''''}) = s(L_{2}^{''''}) = s(L_{3}^{''''})$. Hence by lemma 1.1, $G$ is $S$-magic. Since the highest label is $n + 2 = n + q + 1$, the magic index, $\theta(G) = 2 = q+1$. Since $\lambda - 1 = n_{1} + r - 1$, $\theta(G) = \lceil \frac{\lambda - 1}{n_{1}} \rceil$.\\
    
    \noindent Now, take $q = 0$. Then $\lambda = r$. Since $\lambda > 0$, we have $r > 0$. If $r = 1$, we replace $n \in L_{1}^{'}$ with $n + 1$. Hence the resulting labels are $L_{1}^{v} = \{n - n_{1}, n - n_{1} + 2,\cdots, n - 1, n + 1\}$, $L_{2}^{v} = L_{2}$ and $L_{3}^{v} = L_{3}$. Therefore, $s(L_{1}^{v}) = s(L_{1}) + 1$, $s(L_{2}^{v}) = s(L_{2})$ and $s(L_{3}^{v}) = s(L_{3})$. We label the partitions $V_{1}$, $V_{2}$ and $V_{3}$ with the label sets $L_{1}^{v}$, $L_{1}^{v}$ and $L_{1}^{v}$, respectively. From \eqref{8} and \eqref{9}, it is clear that $s(L_{1}^{v}) = s(L_{3}^{v}) = S(L_{2}^{v})$. Hence by Lemma 1.1, $G$ is $S$-magic. Since the highest label is $n + 1$, the magic index, $\theta(G) = 1 = q + 1$. Since $\lambda - 1 = r - 1$, $\theta(G) = \lceil \frac{\lambda - 1}{n_{1}} \rceil$.
    even
    If $r > 1$. we replace the label $n - r + 1 \in L_{1}^{'}$ with $n - r + 1 + r = n + 1$. Hence the resulting labels are $L_{1}^{v} = \{n - n_{1}, n - n_{1} + 2,\cdots,n - r, n - r + 2, \cdots, n + 1\}$, $L_{2}^{v} = L_{2}$ and $L_{3}^{v} = L_{3}$. Therefore, $s(L_{1}^{v}) = s(L_{1}) + r$, $s(L_{2}^{v}) = s(L_{2})$ and $s(L_{3}^{v}) = s(L_{3})$. We label the partitions $V_{1}$, $V_{2}$ and $V_{3}$ with the label sets $L_{1}^{v}$, $L_{1}^{v}$ and $L_{1}^{v}$, respectively. From \eqref{8} and \eqref{9}, one can see that  $s(L_{1}^{v}) = s(L_{3}^{v}) = s(L_{2}^{v})$. Hence by Lemma 1.1, $G$ is $S$-magic. Since the highest label is $n + 1$, the magic index, $\theta(G) = 1 = q + 1$. Since $\lambda - 1 = r - 1$, we have $\theta(G) = \bigl\lceil \frac{\lambda - 1}{n_{1}} \bigr\rceil$.
    \begin{equation*}
        \therefore \theta(G) = \biggl\lceil \frac{6h - n(n+1) - 2}{2n_{1}} \biggr\rceil = \biggl\lceil \frac{3(n - n_{1})(n - n_{1} + 1) - 2n(n+1) + 2}{4n_{1}} \biggr\rceil
    .\end{equation*}
    Since the subgraph $H$ is $S$-magic with $\theta(H) = 1$, by Corollary 2.7.1, $h$ is the minimum magic constant for the subgraph $H$. Therefore, $\theta(G)$ is minimum as it is a linear function of $h$.\\
  \noindent \textbf{Proof of case} (iii). \\   Assume that $3\zeta_{n-n_1+1}^{n} < \zeta_1^{n}  < 3\zeta_{1}^{n_3}$. Then
    \begin{equation*}
        \frac{n_{1}(2n - n_{1} + 1)}{2} < \frac{n_{3}(n_{3} + 1)}{2}
    .\end{equation*}
    Hence, from equation (4), $g(0) < 0$. Therefore, by Lemma 1.2,
    \begin{equation}
        \theta(G) \ge \biggl\lceil \frac{|n_{1}(2n - n_{1} + 1) - n_{3}(n_{3}+1)|}{2(n - n_{3})} \biggr\rceil.
    \end{equation}
    
\begin{example} Consider the complete tripartite graph $K_{5,6,7}$. Now, by the case $(i)$ of Theorem 2.2. its distance magic index $0$ with magic constant $k = 114$ and it is illustrated in Figure 1. 
    \begin{figure}[ht]
        \centering
        \includegraphics[scale=0.1]{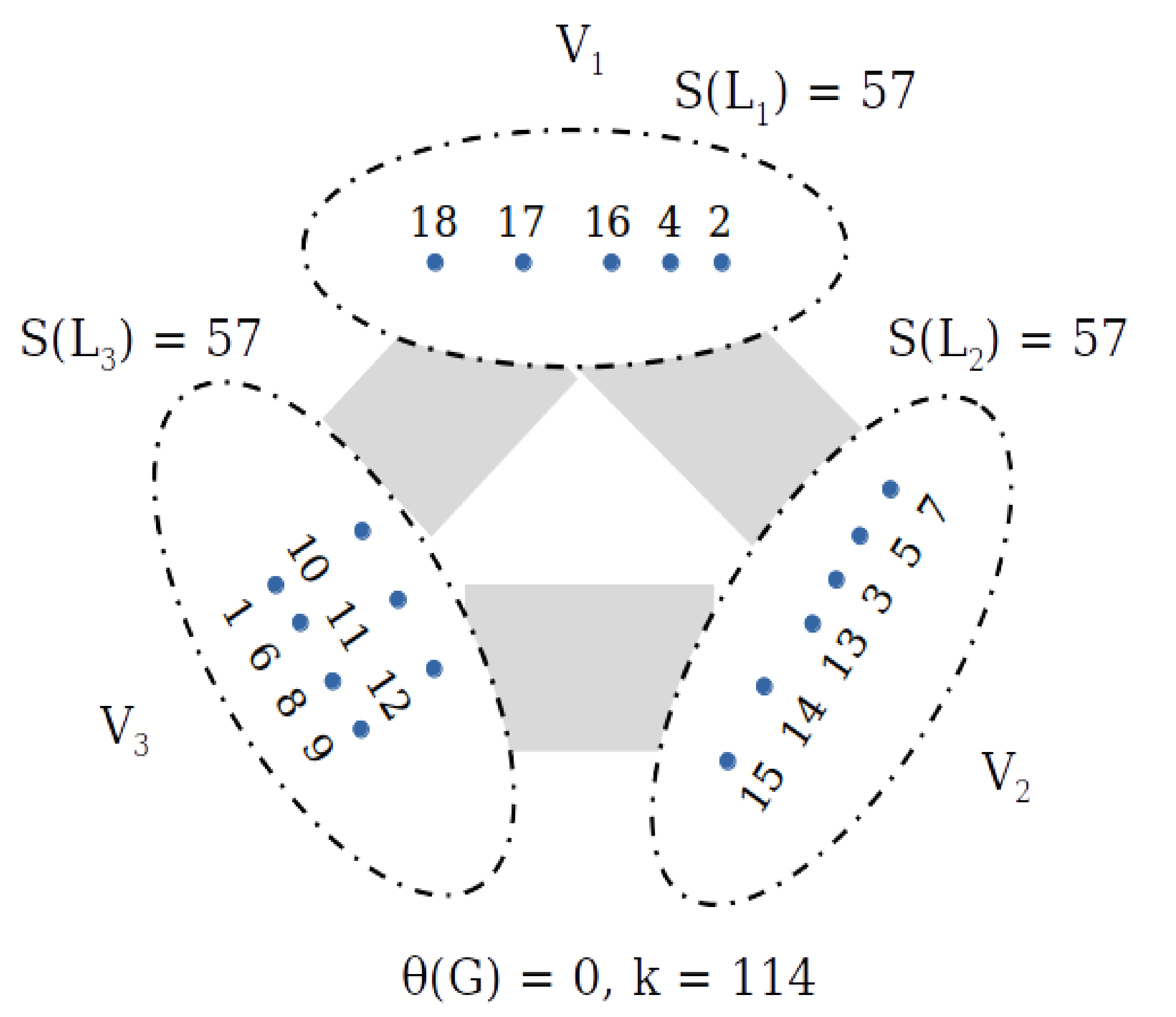}
        \caption{Example of case $(i)$, $G \cong K_{5, 6, 7}$, $\theta(G) = 0$, $k = 114$}
    \end{figure}
   On the other hand, Figure 2 gives an $S$-magic labeling of the complete tripartite graph $K_{3,8,9}$ and its distance magic index $7$ with magic constant $k = 154$. Note that by using case $(ii)$ of Theorem 2.2, one can say that the magic index of $K_{3, 8, 9}$ cannot be $0, 1, 2, 3, 4, 5$ and $6$,
       \begin{figure}[ht]
        \centering
    \includegraphics[scale=0.1]{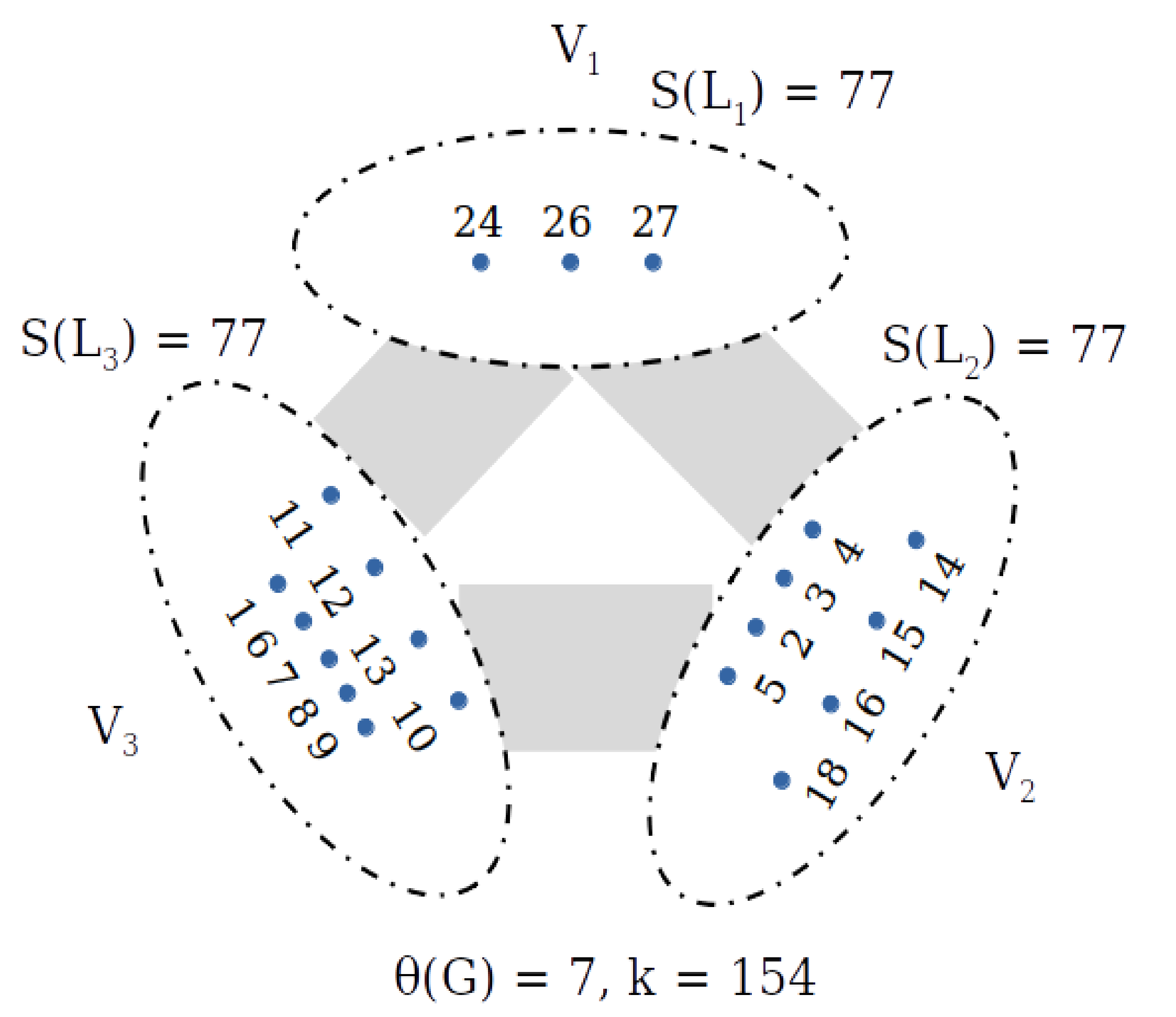}
        \caption{Example of case $(ii)$, $G \cong K_{3, 8, 9}$, $\theta(G) = 7$, $k = 154$}
    \end{figure}\\
    \end{example}
    \noindent \textbf{Proof of case $(iv)$}\\
    Here, $3\zeta_{n-n_1+1}^{n} \ge  \zeta_1^{n} \ge 3\zeta_{1}^{n_3}  \text{ and } 2\zeta_1^{n}  \equiv 2 \text{ mod } 6$. Hence by Theorem 1.3, $\theta(G) \ge 1$. By  Lemma 3.1, one can see that $\theta(K_{n_1 + 1, n_2, n_3}) = 0$.
    Thus, by removing the vertex with maximum label, say $u$, from $V_1$ and adding that label to any other vertex, say $u'$, of $V_1$ such that $u + u' > n$ will make $K_{n_1, n_2, n_3}$, $S$-magic. Since $u \le n + 1$, and $u' \le n$, the highest label, $u + u' \le 2n +1$. Therefore, $\theta(G) \le 2n + 1 - n = n +1$. Hence,
    $1 \le \theta(G) \le n + 1.$
    \newline \noindent But we have to prove the existence of vertices $u$ and $u'$ in $V_1$ such that $u + u' > n$.
    Suppose $u + u' \le n$. Let us consider $S(L_1)$ of graph $K_{n_1+1,n_2,n_3}$
        \begin{eqnarray*}
        s(L_1) &<& (n_1)u + u' \le (n_1-1)u + n\\
            &<& (n_1-1)\frac{n}{2} + n = 
            \frac{n(n_1 + 1)}{2}\\
            &\le& \frac{n(n+3)}{6} < \frac{(n+1)(n+2)}{6}.
        \end{eqnarray*}
        Which is a contradiction, Since $s(L_1) = \frac{(n+1)(n+2)}{6}$.\\
    \textbf{Proof of case $(v)$}.\\
    In this case, we have, $\zeta_1^{n} \le 3\zeta_{n-n_1+1}^{n}$ and $\zeta_1^{n} < 3\zeta_{1}^{n_3} $. Suppose that the partitions $V_1$, $V_2$ and $V_3$ are labelled by the set of labels $L_{1} = \{n,\cdots, n - n_{1} + 1\}$, $L_{2} = \{n - n_{1},\cdots, n_{3} + 1\}$ and $L_{3} = \{n_{3},\cdots, 1\}$, respectively. Let $H \cong K_{n_{2}, n_{3}}$. Consider,
      \begin{eqnarray*}
        \frac{(n-n_{1})(n - n_{1} + 1)}{2} &=& \frac{n^2 + n}{2} - \frac{n_{1}(2n - n_{1} + 1)}{2}\\
        &\le& \frac{n^2 + n}{3} \\ &<& n_{3}(n_{3} + 1).
    \end{eqnarray*}   
    Hence, by using Theorem 1.4, we have
    \begin{equation}
        \theta(H) = \biggl\lceil \frac{|(n - n_{1})(n - n_{1} + 1) - 2n_{3}(n_{3}+1)|}{2n_{2}} \biggr\rceil
    .\end{equation}
            Let $$\frac{|(n - n_{1})(n - n_{1} + 1) - 2n_{3}(n_{3}+1)|}{2} = n_{2}\mu + r, \quad 0 \le r < n_{2}$$
            Hence $H$ is $S$-magic with $L_{3}^{'} = L_{3}$ and\\
            $L_{2}^{'} =$ $\begin{cases}\{n_{3} + 1 + \mu, \cdots, n - n_{1} + \mu\}, & \text{ if } r = 0 \\
            \{n_{3} + 1 + \mu, \cdots, n - n_{1} + \mu - r, n - n_{1} + \mu - r + 2,\cdots, n - n_{1} + \mu + 1\}, & \text{ if } r > 0.\end{cases}$
            Let $L_{1} ^{'} = \{n - n_{1} + \mu + 2, \cdots, n + \mu + 1\}$. Implies, $s(L_{1}^{'}) = s(L_{1}) + n_even{1}(\mu+1)$.\\
           {\bf Subcase 1}.  Suppose $s(L_{1}^{'}) \le s(L_{3})$. Let $s(L_{3}) - s(L_{1}^{'}) = n_{1}\lambda + r_{1}$ where, $0 \le r_{1} < n_{1}$. This implies, $s(L_{3}) - s(L_{1}) = n_{1}\lambda + n_{1}(\mu+1) + r_{1}$. Now, increase each label of $L_{1}^{'}$ by $\lambda$ and hence, we get $L_{1}^{''} = \{n - n_{1} + \mu + \lambda + 2,\cdots, n + \mu + \lambda + 1\}$. If $r_{1} = 0$, we can label the vertices of the partition $V_{1}$, $V_{2}$ and $V_{3}$ by the label set $L_{1}^{''}$, $L_{2}^{'}$ and $L_{3}$, respectively. Since $s(L_{1}^{''}) = s(L_{2}^{'}) = s(L_{3})$, by Lemma 1.1, the graph $G$ is $S$-magic. Since the highest label is $n + \mu + \lambda + 1$, the magic index $\theta(G) = \mu + \lambda + 1$. If $r_{1} > 0$, replace $n + \mu + \lambda + 2 - r_{1} \in L_{1}^{''}$, with $n + \mu + \lambda + 2$. Thus the resulting label set $L_{1}^{'''} = \{n - n_{1} + \mu + \lambda + 2, \cdots, n + \mu + \lambda - r_{1}+1, n + \mu + \lambda + 3 - r_{1}, \cdots, n + \mu + \lambda + 2\}$. Hence, $s(L_{1}^{'''}) = s(L_{1}^{'}) + n_{1}\lambda + r_{1} = s(L_{2}^{'}) = s(L_{3})$. Therefore, by Lemma 1.1, we label the partitions $V_{1}$, $V_{2}$ and $V_{3}$ respectively by $L_{1}^{'''}$, $L_{2}^{'}$ and $L_{3}$ in such a way that the graph $G$ is an $S$-magic. Since the highest label is $n + \mu + \lambda + 2$, the magic index, $\theta(G) = \mu + \lambda + 2$. Therefore, 
    \begin{equation*}
        \theta(G) = \biggl\lceil \frac{2h - n_{1}(2n-n_1+1)}{2n_{1}} \biggr\rceil = \biggl\lceil \frac{n_{3}(n_{3} + 1) - n_1(2n-n_1 + 1)}{2n_{1}} \biggr\rceil
    .\end{equation*}
   {\bf Subcase 2}. Suppose $s(L_{1}^{'}) > s(L_3)$. Hence the maximum highest possible label, is $n + \mu + 1$. Therefore, $\theta(G) \le \mu + 1 \le \theta(H) + 1$ and hence $\theta(G) \le \theta(H) + 1$.
    \end{proof}

\par  Now, we discuss the distance magic index of complete multi-partite graphs and some related graphs. We denote $K(a,b)$ as the complete multi-partite graph with $b$ partitions and each partition having $a$ vertices in it, where $a\geq 1, b\geq 1$. Note that for a graph $G$, the graph $mG$ is the disjoint union of $m$ copies of the graph $G$. Recall that the lexicographic product~\cite{MR2817074} of two graphs $G$ and $H$ is a graph $G \circ H$ with the vertex set $V(G)\times V(H)$ and two vertices $(u,v)$ and $(u',v')$ are adjacent if and only if $u$ is adjacent to $u'$ in $G$, or $u=u'$ and $v$ is adjacent to $v'$ in $H$.
\begin{theorem}
	Let $h \geq 1.$ A Kotzig array $KA(2h, b)$ exists for any $b$;
	$KA(2h + 1, b)$ exists if and only if $b$ is odd.
\end{theorem}
In 2022, Froncek et al.\cite{AVPrajeesh}  introduced the notion of the lifted quasi-Kotzig arrays and used them to construct the quasimagic rectangles as a generalization of magic rectangles. 
	\begin{definition} \cite{AVPrajeesh}
		For an odd $a$ and even $b$, a quasi-magic rectangle $QMR(a,b:d)$ is an $a\times b$ array with the entries $1,2,\dots,d-1,d+1,\dots, ab+1$, each appearing exactly once, such that the sum of each row is equal to a constant $\rho$ and the sum of each column is equal to a constant $\sigma$.\end{definition}
	They proved the existence of $QMR(a,b:d)$ for all possible values of $a$ and $b$ and $d=ab/2+1$. When $\gcd(a,b)=1$, and $d=ab/2+1$, they obtained a characterization for the existence of $QMR(a,b:d)$.
	\begin{theorem}\textnormal{\cite{AVPrajeesh}}\label{thm:main}
		{There exist quasimagic rectangles $QMR(a,2t:a t+1)$ for all odd $a$ and $t \geq 1$ with exceptions of $QMR(a,2:a+1)$ for $a \equiv 1(\bmod~4)$. In addition, when $\gcd(a, 2 t)=1$,  then the necessary condition is also sufficient.}
	\end{theorem}
	Note that it is certainly, one of the tools to determine the distance magic index of certain partite graphs.

\section{Quasimagic rectangle: A tool to determine the distance magic index}
\par  In this section, the distance magic indices of the graphs such as $mK(a,b)$, $m(C_{b}\circ\overline{K}_a)$ and $G\circ \overline{K}_a$, are computed, where $G$ is an arbitrary regular graph. These graphs already appear in the literature in the articles \cite{MR2848089,MR1999203,MR2533199}. These articles mainly discuss about the distance magicness of the above mentioned graphs. The existing results of  distance magicness of these graphs, are tabulated below for a quick reference. \\
\noindent Case 1: The distance magicness of $K(a,b)$ is given in Table \ref{Table1}.
\begin{table}[h]
	\centering
	\begin{tabular}{|c|c|c|}
		\hline
		values of $a$ & odd $b$'s & even $b$'s\\
		\hline
		even & $dmg$ \cite{MR1999203} & $dmg$ \cite{MR1999203} \\
		\hline
		odd & $dmg$ \cite{MR1999203} & not~  $dmg$ \cite{MR1999203}  \\
		\hline
	\end{tabular} 
	\caption{}  
	\label{Table1}
\end{table}
\\ \noindent Case 2: The distance magicness of $mK(a,b),m>1$ is given in Table \ref{Table2}.
\begin{table}[htbp]
 \centering
\begin{tabular}{|c|c|c|c|}
\hline values of $m$ & values of $a$ & odd $b$'s & even $b$'s\\ \hline
odd  & even & $dmg$ \cite{MR2533199}  &$dmg$\cite{MR2533199}\\ \hline
even & even  &  $dmg$ \cite{MR2533199}& $dmg$ \cite{MR2533199} \\	\hline
odd &  odd  &$dmg$ \cite{MR2533199}  &  not $dmg$ \cite{MR2533199}  \\ \hline
 even & odd &  not dmg if $b=4t+1$ \cite{MR2848089}, not dmg if $b=4t+3$ \cite{MR2533199} 
& not $dmg$ \cite{MR2533199} \\ \hline 
\end{tabular}
\caption{}  
\label{Table2}
\end{table}
\newline
Case 3: The distance magicness of $m(C_b\circ\overline{K}_a)$ is given in Table \ref{Table3}.
\begin{table}[htbp]
	\begin{center}
		{\footnotesize \begin{tabular}{|c|c|c|c|}
				\hline
				values of $m$ & values of $a$ & even $b$'s & odd $b$'s\\
				\hline
				odd &  {even} &  $dmg$ \cite{MR2533199}& $dmg$ \cite{MR2533199} \\\hline
				even & even &  $dmg$ \cite{MR2533199} & $dmg$ \cite{MR2533199} \\
				\hline
				odd &  odd   & $dmg$ if $b=4t$, {not $dmg$} if $b=4t+2$  \cite{MR2533199}&  {$dmg$} \cite{MR2533199} \\ \hline
				even & odd & $dmg$ if $b=4t$, {not $dmg$} if $b=4t+2$ \cite{MR2533199}  &  {not $dmg$} \cite{MR2533199} \\
				\hline
		\end{tabular}}
	\end{center}
	\caption{}
	\label{Table3}
\end{table}
\newline
Case 4: If $G$ is an $r$-regular graph on $b$ vertices, the distance magicness of $G\circ\overline{K}_a$ is given in Table \ref{Table4}.
\begin{table}[h]
	\begin{center}
		{\footnotesize \begin{tabular}{ |c|c|c|c|c|}
				\hline
				values of $r$
				& values of $a$ & even $b$'s, $b=4t$ & even $b$'s, $b=4t+2$& odd $b$'s\\
				\hline
				even, $r=4s$ &  odd   & {not yet solved}& { not yet solved} & $dmg$ \cite{MR2848089}\\\hline
				even, $r=4s+2$ & odd & { not yet solved } &  {not $dmg$} \cite{MR2848089} &  $dmg$ \cite{MR2848089}\\hline
				odd & & {not $dmg$}& {not $dmg$}& no such graphs\\\hline
				even & even  & $dmg$\cite{MR1999203} &  $dmg$ \cite{MR1999203}
				&  $dmg$ \cite{MR1999203}
				\\\hline
				odd &even&  $dmg$ \cite{MR1999203}& $dmg$ \cite{MR1999203}& no such graphs \\ 
				\hline
		\end{tabular}}
	\end{center}
	\caption{} 
	\label{Table4}
\end{table}
\begin{theorem}\label{index1hnp} For the graph $K(a,b)$, 	\\
	$\theta(K(a,b)) =$	
	$\begin{cases}
	0 &   \text{for $a$ even 
		or $a$ and $b$ both odd} \\
	1 &   \text{for $a$ odd and $b$ even except when $a \equiv 1 \bmod 4$ and $b=2$}.
	\end{cases}$
\end{theorem} 
\begin{proof}
	If $a$ is even or $a$ and $b$ are both odd, then from Table 1, $\theta(K(a,b))=0$. Now, if  $a$ is odd and $b$ is even except when $a \equiv 1 \mod 4$ and $b=2$, then again from Table 1,   $\theta(K(a,b))\neq 0$.
	\par Now, if $a$ is odd and $b$ is even, then by using Theorem 2.4, construct quasimagic rectangle $QMR(a,b)$ from  the set $S = \{1,2,\cdots,\frac{ab}{2},\frac{ab}{2}+2,\cdots,ab+1\}$. Then, label $i^{th}$  partition of $K(a,b)$ by $i^{\textnormal{th}}$ column of $QMR(a,b)$, where $1\leq i \leq b$. Thus, if $\rho$ is the column sum of a $QMR(a,b)$, then $K(a,b)$ is $S$-magic graph with magic constant $\rho(b-1)$ and therefore, $\theta(K(a,b))$ is 1.
\end{proof}
\begin{theorem}\label{multi6}
	Let $a>1,b>1,m> 1$. Then\\
	$\theta(mK(a,b))$ = 
	$\begin{cases}
	0 &   \text{for} a \text{is even $(or)$} 
              ~mab \text{~is odd}\\
	1 &  \text{otherwise.}
	\end{cases}$ 
\end{theorem}
\begin{proof}
	From Table 2, it can be observed that when $m$ is even, and $a$ is odd or $m,a$ are odd and $b$ is even; we have $\theta(mK(a,b))\neq 0$. For the remaining cases $\theta(mK(a,b))= 0$.
	\par For all the cases, when $mK(a,b)$ is not distance magic, by using Theorem 2.4,  construct a $QMR(a,mb)$ with column sums $\rho$.
	Now, use it to label $m$ copies of the graph $K(a,b)$ by choosing $mb$ columns of $QMR(a,mb)$ to the $(mb)$-partitions of $mK(a,b)$.  Hence, an $S$-magic labeling of $mK(a,b)$ with $\mu' =\rho(b-1)$ is obtained. Therefore, $\theta(mK(a,b))=1.$
\end{proof}
\begin{theorem}\label{thm14}
	Let $m\geq 1, a>1$ and $b\geq3$.  Then\\
	\begin{equation*}
	\theta(m(C_{b}\circ \overline{K}_a) = 
	\begin{cases}
	0 &   \text{if $a$ is~even $(or)$ $mab$ is odd, $(or)$ $a$ is odd, and $b\equiv 0 \textnormal{~mod~} 4$},\\
	1 &   otherwise.
	\end{cases}  
	\end{equation*}
\end{theorem}
\begin{proof}
	If $a$ is even or $mab$ is odd, or $a$ is odd and $b \equiv 0 \bmod 4$, then from Table 3, $\theta(G)=0.$
	Now the remaining cases are given below. \\
	Case 1. $a$ is odd, $m$ is even, $b\equiv 2 \bmod 4.$\\
	Case 2. $a$ is odd, $m$ is odd, $b\equiv 2 \bmod 4.$\\
	Case 3. $a$ is odd, $m$ is even, $b$ is odd.
	\par Now, for all these cases, a $QMR(a,mb)$ is constructed by using Theorem 2.4, with individual column-sum $\rho$ and individual row-sum $\sigma$. The first $b$ columns of $QMR(a,mb)$ is used to label the first copy of  $m(C_{b}\circ \overline{K}_a)$, by labeling   the $i^{th}$ layer of $K_{a}$ of $C_{b}\circ \overline{K}_a$ using the $i^{th}$ column of $QMR(a,mb)$. Similarly, the next $b$ columns are used to label the vertices of the second copy  $C_{b}\circ \overline{K}_a$. If we proceed like this, one can label all copies of $m(C_{b}\circ \overline{K}_a)$, using the column of  $QMR(a,mb)$.  Thus, $m(C_{b}\circ \overline{K}_a)$ is $S$-magic with $\mu' = 2\rho$ and therefore, $\theta(m(C_{b}\circ \overline{K}_a))=1$.
\end{proof}
\begin{theorem}\label{thm15}
	Let $a\ge 1$, and $G$ be an $r$-regular graph on $b$ vertices. Then
	\begin{equation*}
	\theta(G\circ \overline{K}_a) =
	\begin{cases}
	1 &   \text{for $a$ and $r$ are odd except when $a \equiv 1 \bmod 4$ and $b=2$} \\
	1 &   \text{for $a$ is odd, $r\equiv 2 \bmod 4$, $b \equiv 2\bmod 4$},\\
	0 &   \text{for $a$ is even $(or)$ $a$ and $b$ are odd, and $r$ is even}.
	\end{cases}
	\end{equation*}
	
\end{theorem}
\begin{proof}
	Let $v_{1},v_2,\cdots,v_{b}$ be the vertices of $G$. Now, for any $j\in\{1,2,\cdots,b\}$, let $V_{j} = \{v_{j}^{i}: i\in\{1,2,\cdots, a\}\}$ be set the vertices of $G\circ \overline{K}_a$ that replace the vertex $v_{j}$ of $G$. In fact, the vertex set of $G\circ \overline{K}_a$ is $\displaystyle \bigcup_{j=1}^{b}V_{j}$ and $G\circ \overline{K}_a$ is a $(ra)$-regular graph on $ab$ vertices. \par Now consider the following cases:\\ Case 1. $a$ is even. From Table 4, we get $\theta(G\circ \overline{K}_a)=0$.\\
	Case 2. $a$ is odd. When $r$ is even, and $b$ is odd, from Table 4, we have
	$\theta(G\circ \overline{K}_a)=0$. On the other hand, if $r\equiv 2 \bmod 4$, and $b \equiv 2 \bmod 4$, then from Table 4,   $\theta(G\circ \overline{K}_a)\neq 0$. Finally, if $r$ is odd, then $b$ is even, and by Theorem \ref{oddregular}, $\theta(G\circ \overline{K}_a)\neq0$.\\
	\par Now, for both the cases, when $\theta(G\circ \overline{K}_a)\neq0$,  construct $QMR(a,b)$ using Theorem 2.4. Use the $j^{th}$ column of $QMR(a,b)$ to label the set of vertices, $V_{j}$, for all $j \in \{1,2,\cdots,b\}$. Hence, an $S$-magic labeling of $G\circ \overline{K}_a$, with  $\mu' = r\rho$ is obtained. Therefore, $\theta(G\circ \overline{K}_a)=1.$
\end{proof} 
\begin{example}
		A quasimagic rectangle $QMR(3,10)$ is given below. Now, Figure \ref{Fig5.3} illustrates how one can use this $QMR$ to give $S$-magic labeling of the $9$-regular graph $G \circ \overline{K}_3$, where $S$ is the set of all entries of this $QMR$.  \vskip 0.12cm
\begin{center}
\begin{tabular}{|c|c|c|c|c|c|c|c|c|c|}
		\hline	
   22 & 23 & 5 &  20 & 1 & 10 &  12 & 6 & 30 & 31 \\\hline	
		17 & 18 & 19 & 3 & 26 & 27 &  8 & 13 & 14 & 15 \\\hline	
		9 & 7 & 24 & 25 & 21 & 11 & 28 & 29 & 4 & 2  \\\hline
\end{tabular}
 \end{center}
\vskip .2cm \noindent One can see that the column sums are 48 and row sums are 160 and moreover, $\theta(G\circ \overline{K}_3)=1.$
\begin{figure}[ht]
	\centering
	\includegraphics[width=120mm]{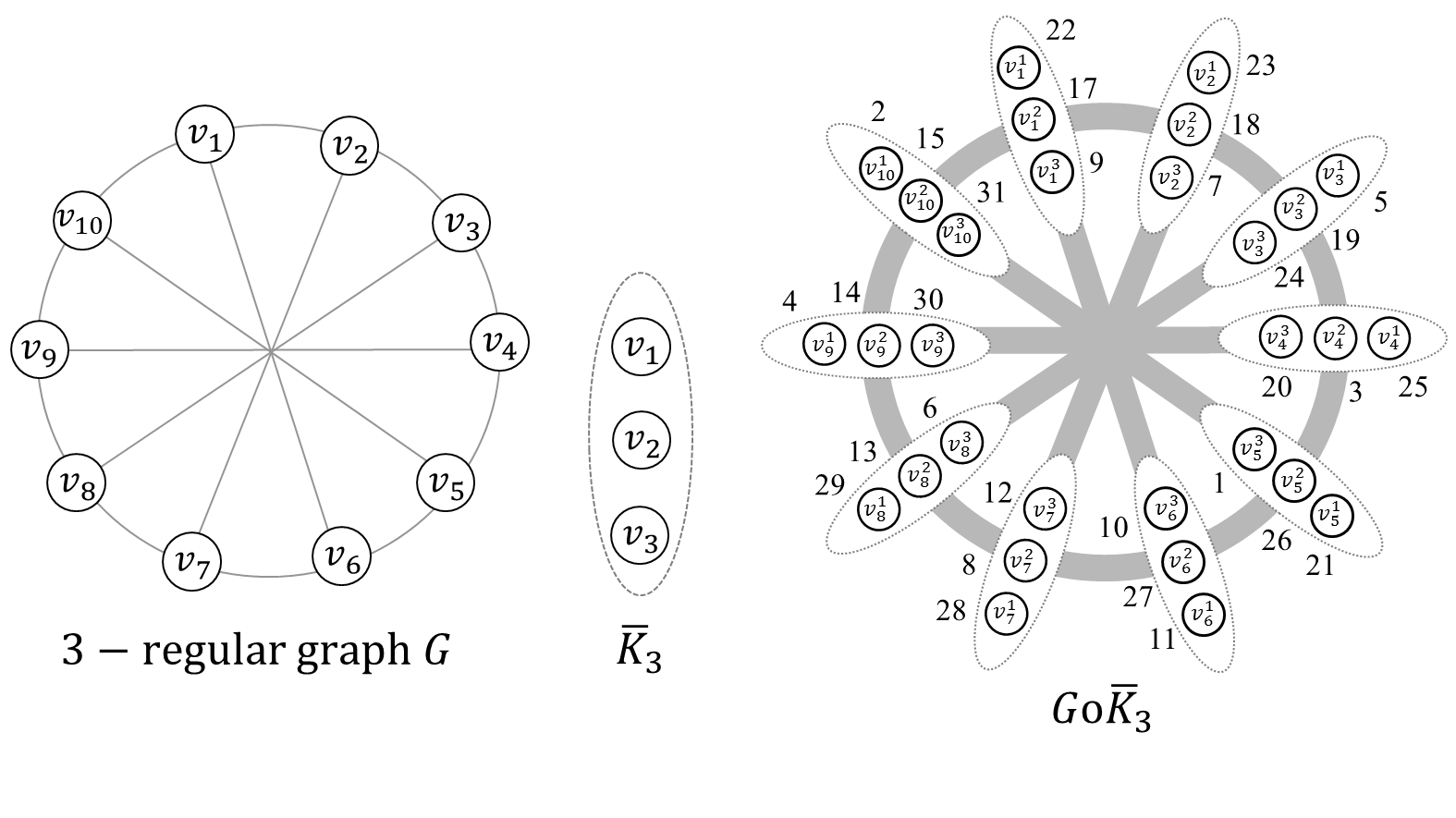}
	\caption{An $S$-magic labeling of $G \circ \overline{K}_3$ using $QMR(3,10)$}
	\label{Fig5.3}
\end{figure}
 \end{example}

\section{Conclusion and scope}
In this paper,   the distance magic index of complete tri-partite graphs $K_{n_1,n_2,n_3}$ for all possible values of $n_1,n_2$ and $n_3$, are determined and  Table 5 gives the summary. However, the problem of determining the distance magic index of complete $r$-partite graphs $K_{n_1,n_2,\cdots, n_r}$ is still open, where $r>3$. In last section, the quasimagic rectangle $QMR$ is  identified as useful tool to determine the distance magic indices of certain regular graphs such as disjoint union of $m$ copies of $K(a,b)$, disjoint union of $m$ copies of $C_{b}\circ \overline{K}_a$ and $G \circ {K}_a$, where $G$ is an arbitrary regular graph.


\begin{thebibliography}{}

\bibitem{MR3082443} S. Arumugam,  D. Froncek, and N. Kamatchi, Distance magic graphs--a survey,
Journal of the Indonesian Mathematical Society, Special edition, 11-26, 2011.

\bibitem{MR2368647} J.A. Bondy, and U.S.R. Murty, Graph Theory, Springer, New York, 2008.  


\bibitem{AVPrajeesh} D. Fron\v{c}ek,  K. Paramasivam, K. and A V Prajeesh,  Quasimagic rectangles, Journal of Combinatorial Designs, 30(3), 193-202, 2022.

\bibitem{MR2848089} D. Fron\v{c}ek,  P. Kov\'{a}\v{r}, and T. Kov\'{a}\v{r}ov\'{a}, Constructing distance magic graphs from regular graphs, Journal of Combinatorial Mathematics and Combinatorial Computing, 78, 349-354, 2011.

\bibitem{MR2259699} D. Fron\v{c}ek,  P. Kov\'{a}\v{r}, and T. Kov\'{a}\v{r}ov\'{a}, Fair incomplete tournaments, Bulletin of the Institute of Combinatorics and its Applications, 48, 31-33, 2006
\bibitem{MR3803223} A. Godinho, and T. Singh, Some distance magic graphs, AKCE International Journal of Graphs and Combinatorics, 15(1), 1-6, 2018
    
\bibitem{MR3743955} A. Godinho, T. Singh, and S. Arumugam, The distance magic index of a graph, Discussiones Mathematicae Graph Theory, 38(1), 135--142, 2018,

 

\bibitem{godinho2015s}  A. Godinho,  T. Singh, and S. Arumugam, On $S$-magic graphs, Electronic Notes in Discrete Mathematics, 48,267-273,2015.

 
 \bibitem{MR1710483} T. R. Hagedorn,  Magic rectangles revisited, Discrete Mathematics, 207, 65-72, 1999.
 
\bibitem{MR2817074} R. Hammack, W. Imrich, and S. Klav\v{z}ar, Handbook of Product Graphs,CRC Press, Boca Raton, FL, 2011.


\bibitem{kotzig1971magic} A. Kotzig, On magic valuations of trichromatic graphs. Reports of the CRM, 1971.

\bibitem{MR3013245} P. Kov\'{a}\v{r}, D. Fron\v{c}ek, and T. Kov\'{a}\v{r}ov\'{a}, A note on 4-regular distance magic graphs, The Australasian Journal of Combinatorics, 54, 127--132, 2012,

\bibitem{kotlar} D. Kotlar, Distance magic labeling in complete $4$-partite graphs, Graphs and Combinatorics,32(3), 1027-1038, 2016.

\bibitem{MR3013201} A. M. Marr, and W.D. Wallis, Magic Graphs, Birkh\"{a}user/Springer, New York, 2013.

\bibitem{MR1999203} M. Miller, C. Rodger, and R. Simanjuntak, Distance magic labelings of graphs,
The Australasian Journal of Combinatorics, 28,305--315, 2003.
\bibitem{Rao2004} S. B. Rao, T. Singh, and V. Parmeswaran, Some sigma labelled graphs I, Proc. Graphs, Combinatorics, Algorithms and Applications, 125-133, 2004.
\bibitem{MR2533199} M. K. Shafiq, G. Ali, and R. Simanjuntak, Distance magic labelings of a union of graphs, AKCE International Journal of Graphs and Combinatorics, 6(1), 191--200, 2009.

\bibitem{vilfred} V. Vilfred,  $\sum$-labelled graphs and circulant graphs,  Ph.D. Thesis, University of Kerala,  India, 1994.
\end{thebibliography}

\end{document}